\documentclass{article}
\begin{document}
\centerline{\textbf{\Large{Eisenstein Series, Alternative Modular Bases}}}
\centerline{\textbf{\Large{and}}}
\centerline{\textbf{\Large{Approximations of $1/\pi$}}}
$$
$$
\centerline{\bf Nikos Bagis}

\centerline{Department of Informatics}

\centerline{Aristotle University of Thessaloniki}
\centerline{54124 Thessaloniki, Greece}
\[
\]
\begin{quote}
\begin{abstract}
In this article using the theory of Eisenstein series, we give rise to the complete evaluation of two Gauss hypergeometric functions. Moreover we evaluate the modulus of each of these functions and the values of the functions in terms of the complete elliptic integral of the first kind. As application we give way of how to evaluate the parameters, in a  closed-well posed form, of a general Ramanujan type $1/\pi$ formula. The result is a formula of 110 digits per term. 
\end{abstract}

\bf keywords \rm{elliptic functions; singular modulus; Ramanujan's cubic theory; Pi}

\end{quote}

\section{Introduction}
\label{intro}
1) The Gauss hypergeometric function is (see [1]) 
\begin{equation}
_2F_1(a,b;c;w):=\sum^{\infty}_{n=0}\frac{(a)_n(b)_n}{(c)_n}\frac{w^n}{n!}
\end{equation} 
where $(a)_n=\frac{\Gamma(n+a)}{\Gamma(a)}=a(a+1)(a+2)\ldots(a+n-1)$, and $|w|<1$.\\
2) The Legendre function (see [1]) $$P^{\mu}_{\nu}(w)=\frac{1}{\Gamma(1-\mu)}\left(\frac{w+1}{w-1}\right)^{\mu/2}{}_2F_1\left(-\nu,\nu+1;1-\mu;\frac{1-w}{2}\right)$$ is solution of the equation $$(1-w^2)\frac{d^2y}{dw^2}-2w\frac{dy}{dw}+\left(\nu(\nu-1)-\frac{\mu^2}{1-w^2}\right)y=0$$
3) Let also $q=e^{-\pi\sqrt{r}}$, where $r$ positive real. 
The cubic theta functions are (see [6]) 
\begin{equation}
a(q):=\sum^{\infty}_{m,n=-\infty}q^{m^2+mn+n^2}     
\end{equation}
\begin{equation}
c(q):=\sum^{\infty}_{m,n=-\infty}q^{(m+1/3)^2+(m+1/3)(n+1/3)+(n+1/3)^2}
\end{equation}
and the cubic singular modulus is given by
\begin{equation}
\alpha_r=\left(\frac{c(q^{\frac{2}{\sqrt{3}}})}{a(q^{\frac{2}{\sqrt{3}}})}\right)^3
\end{equation}
Then 
\begin{equation}
\frac{_2F_1\left(\frac{1}{3},\frac{2}{3};1;1-\alpha_r\right)}{_2F_1(\frac{1}{3},\frac{2}{3};1;a_r)}=\sqrt{r}
\end{equation}
\begin{equation} 
P_{-1/3}(1-2w)={}_2F_1(1/3,2/3;1;w)=z
\end{equation}
The function 
\begin{equation}
z:={}_2F_1\left(\frac{1}{3},\frac{2}{3};1;w\right)=\int^{1}_{0}\frac{1}{\sqrt[3]{t(1-t)^2(1-tw)}}dt
\end{equation}
is called cubic elliptic function and is similar in properties to those of the elliptic integral of the first kind $K(w)$ (see [1],[4],[8],[9],[10]):
$$ _2F_1\left(\frac{1}{2},\frac{1}{2};1;w\right)=\frac{2}{\pi}\int^{\pi/2}_{0}\frac{1}{\sqrt{1-w \sin^2(t)}}dt=\frac{2}{\pi}K(w)
$$
We set $m_r$ to be the solution of $$\frac{K(1-m_r)}{K(m_r)}=\sqrt{r}$$
and call the $k_r=(m_r)^{1/2}$ singular moduli.\\
As with the elliptic integral $K$, using the cubic elliptic theory we give evaluations of the singular modulus $\alpha_r$ for certain $r$ (see [5],[6]). Also we derive a formula for finding  $$_2F_1\left(\frac{1}{3},\frac{2}{3};1;a_r\right)$$
4) We also consider and treat with  $u(w):=P_{-1/6}(w)={}_2F_1\left(\frac{1}{6},\frac{5}{6};1;w\right)$ as in the cubic case. The $\beta_r$ are solutions of the equation $$\frac{_2F_1\left(\frac{1}{6},\frac{5}{6};1;1-\beta_r\right)}{_2F_1\left(\frac{1}{6},\frac{5}{6};1;\beta_r\right)}=\sqrt{r}$$ 
The cases 
$_2F_1\left(\frac{1}{2},\frac{1}{2};1;w\right)$ and $_2F_1\left(\frac{1}{4},\frac{3}{4};1;w\right)$ are trivial and easily can studied. Moreover the first is $K$ itself. There are also other hypergeometric functions that can be evaluated with $K$. 
For example
$${}_2F_1\left(\frac{1}{4},\frac{7}{4};1;w\right)=-\frac{2 \left(E\left(2-\frac{2}{1+\sqrt{w}}\right)-2 \left(-1+\sqrt{w}\right) K\left(2-\frac{2}{1+\sqrt{w}}\right)\right)}{3\pi\left(-1+\sqrt{w}\right) \sqrt{1+\sqrt{w}}} \eqno{:(a)}$$ 
Where $E$ is the complete elliptic integral of the second kind.\\  
In most of our results we use Mathematica and in some cases we have no direct proofs.  

\section{Eisenstein Series and the Cubic base} 

If we set $q_1=e^{2\pi i\tau}$, $\tau=\sqrt{-r}$ then (see [5],[6]):   
\begin{equation}
g_2^{*}(\tau)=\frac{4\pi^4}{3}(1+8\alpha)z^4
\end{equation}
\begin{equation}
g_3^{*}(\tau)=\frac{8\pi^6}{27}(1-20\alpha-8\alpha^2)z^6
\end{equation}
where 
\begin{equation}
g^{*}_2(\tau)=60\sum^{\infty}_{(0,0)\neq (m,n)=-\infty}\frac{1}{(m\tau+n)^4}
\end{equation}
\begin{equation}
g_3^{*}(\tau)=140\sum^{\infty}_{(0,0)\neq (m,n)=-\infty}\frac{1}{(m\tau+n)^6}
\end{equation} 
or from the Theory of Eisenstein Series holds
\begin{equation}
g_\nu(r)=2\zeta(2\nu)-\frac{8\nu\zeta(2\nu)}{B_{2\nu}}\sum^{\infty}_{n=1}\frac{n^{2\nu-1}e^{2\pi in\tau}}{1-e^{2\pi i n\tau}}
\end{equation}\
where $\zeta(s)$ is the Riemann's Zeta function and $B_n$ are the Bernoulli numbers, (see [1],[9]).
Hence
$$g_2(r)=\frac{g_2^{*}(\tau)}{60}$$
$$g_3(r)=\frac{g_3^{*}(\tau)}{140}$$
\begin{equation}
\sum^{\infty}_{n=1}\frac{n^{2\nu-1}}{e^{2\pi n\sqrt{r}}-1}=
\frac{2\zeta(2\nu)-g_{\nu}(\sqrt{-r})}{8\nu\zeta(2\nu)}B_{2\nu}
\end{equation}
Hence
\begin{equation}
\sum^{\infty}_{n=1}\frac{n^3}{e^{2\pi n\sqrt{r}}-1}=
\frac{2\zeta(4)-g_2(\sqrt{-r})}{16\zeta(4)}B_{4}\end{equation}
\begin{equation}
\sum^{\infty}_{n=1}\frac{n^5}{e^{2\pi n\sqrt{r}}-1}=
\frac{2\zeta(6)-g_3(\sqrt{-r})}{24\zeta(6)}B_{24}\end{equation}
Now using the above relations we have
\begin{equation} 
\sum^{\infty}_{n=1}\frac{n^3}{e^{2\pi n\sqrt{r/3}}-1}=
-\frac{1}{240}+\left(\frac{1}{240}+\frac{\alpha_r}{30}\right) z^4
\end{equation}
\begin{equation} 
\sum^{\infty}_{n=1}\frac{n^5}{e^{2\pi n\sqrt{r/3}}-1}=
\frac{1}{504}+\left(-\frac{1}{504}+\frac{5\alpha_r}{126}+\frac{\alpha_r^2}{63}\right)z^6
\end{equation}
or
\[
\]
\textbf{Proposition 1.}
\begin{equation} 
1+240\sum^{\infty}_{n=1}\frac{n^3}{e^{2\pi n\sqrt{r/3}}-1}=\left(1+8\alpha_r\right) z^4
\end{equation}
\begin{equation} 
1-504\sum^{\infty}_{n=1}\frac{n^5}{e^{2\pi n\sqrt{r/3}}-1}=
\left(1-20 \alpha_r-8\alpha_r^2\right)z^6
\end{equation}
In the appendix I have construct a table with the Mathematica  of $a_r$ for $r=1,2,\ldots,37$, and $r=40,44,49,59$.\\ 
We set $z_3:={}_2F_1(1/3,2/3;1;\alpha_{3r})$ and we begin to search relations between $F:=\frac{2}{\pi}K$ and $z$, and $u$.
\begin{equation} 
1+240\sum^{\infty}_{n=1}\frac{n^3q^{2n}}{1-q^{2 n}}=\left(1+8\alpha_{3r}\right) z_3^4
\end{equation}
\begin{equation} 
1-504\sum^{\infty}_{n=1}\frac{n^5q^{2n}}{1-q^{2n}}=
\left(1-20 \alpha_{3r}-8\alpha_{3r}^2\right)z_3^6
\end{equation}
It is known that (see [4] pg. 456), (here $x=k_r^2$) 
\begin{equation} 
Q(q^2)=1+240\sum^{\infty}_{n=1}\frac{n^3q^{2n}}{1-q^{2 n}}=(1-x+x^2)F^4
\end{equation}
\begin{equation} 
R(q^2)=1-504\sum^{\infty}_{n=1}\frac{n^5q^{2n}}{1-q^{2n}}=
(1+x)\left(1-\frac{x}{2}\right)(1-2x)F^6
\end{equation}
Also if $y=\pi\sqrt{r}$
\begin{equation}
P(q^2)-\frac{3}{y}=1-\frac{3}{y}-24\sum^{\infty}_{n=1}\frac{nq^{2n}}{1-q^{2n}}=\left(3\frac{E(x)}{K(x)}-2+x-\frac{3 \pi}{4\sqrt{r}K(x)^2}\right)F^2
\end{equation}
Hence the number
\begin{equation} 
s_r=3\frac{E(x)}{K(x)}-2+x-\frac{3 \pi}{4\sqrt{r}K(x)^2}=1+k^2_r-\frac{3a(r)}{\sqrt{r}}
\end{equation}
according to Ramanujan is algebraic, for $r$ positive rational. The function $a(r)$ is called  elliptic alpha function (see [8]). Hence because $x=m_r=k^2_r$ we get the following evaluation of $_2F_1(1/3,2/3;1;\alpha_r)$
\[
\]
\textbf{Theorem 1.}
\begin{equation}
_2F_1\left(\frac{1}{3},\frac{2}{3};1;a_r\right)=\frac{2}{\pi}\sqrt[4]{\frac{1-k_{r/3}^2+k_{r/3}^4}{1+8a_r}}K(k_{r/3})
\end{equation}
\textbf{Proof.}\\
Use (20),(22)
\[
\]
\textbf{Examples.}\\
1)
$$_2F_1\left(\frac{1}{3},\frac{2}{3};1;a_3\right)=8\left(1+\frac{2}{\sqrt{3}}\right)^{1/4}\frac{\sqrt{2\pi}}{\Gamma\left(-\frac{1}{4}\right)^2}$$
2)
$$_2F_1\left(\frac{1}{3},\frac{2}{3};1;a_6\right)=\frac{5}{2\sqrt{6\pi}}\left(\frac{3-2\sqrt{2}}{29-6\sqrt{6}}\right)^{1/4}\frac{\Gamma\left(\frac{1}{8}\right)}{\Gamma\left(\frac{5}{8}\right)}$$
\[
\]
In [4] chapter 21 one can see a very large number of hidden identities very useful for the evaluation of $a(n)$. More precisely using (24) and (25) we can evaluate $a(n^2r)$ in terms of $m_n$, $k_{n^2r}$, $k_r$, $\sqrt{r}$, where $m_n=K(k_{n^2r})/K(k_{r})=K[n^2r]/K[r]$ is the multiplier of the modular equation of $n$-th degree, in the base $K$.
Examples are (see [4] pg. 460 Entry 3.):
$$\frac{a(9r)}{\sqrt{r}}-k^2_{9r}=1-\frac{k_{9r}k_{r}}{3m_3}-\frac{k'_{9r}k'_{r}}{3m_3}-\frac{1}{3m_3}-\frac{1}{3m^2_3}+\frac{1}{m^2_3}\left(\frac{a(r)}{\sqrt{r}}-\frac{k^2_{r}}{3}\right)\eqno{(a1)}$$
$$
\frac{a(81r)}{\sqrt{r}}-3k^2_{81r}=3-\frac{m^5_3}{6m^{7/2}_9}-\frac{m^3_3}{m^{5/2}_9}-\frac{1}{3m^2_9}-\frac{3m_3}{2m^{3/2}_9}+\frac{1}{m^2_9}\left(\frac{a(r)}{\sqrt{r}}-\frac{k^2_r}{3}\right) \eqno{(a2)}
$$
Where $m_3$ is root of the polynomial 
$$27m^4_3-18m^2_3-8(1-2k^2_r)m_3-1=0$$
\[
\]
\textbf{Theorem 2.}\\Given $r$ we set $x=k_r^2$ and $t_1=(1-x+x^2)F^4$, $t_2=(1+x)\left(1-\frac{x}{2}\right)(1-2x)F^6$, $F=\frac{2}{\pi}K(x)$.
The system of equations:
\begin{equation}
(1+8a_{3r})z_3^4=t_1
\end{equation}
\begin{equation}
(1-20a_{3r}-8a_{3r}^2)z_3^6=t_2
\end{equation}
have solutions $(a_{3r},z_{3})$ which can expressed in radicals of $t_1$ and $t_2$. \\  
\textbf{Proof.}\\
Use Mathematica to solve the system with the command 'Solve'. The solution is very complicated to present it here.
\[
\]
\textbf{Theorem 3.}
\begin{equation}
u(\beta_r)={}_2F_1\left(\frac{1}{6},\frac{5}{6};1;\beta_r\right)=\frac{2}{\pi}\sqrt[4]{1-k_{r}^2+k_{r}^4}K(k_{r})
\end{equation}
\textbf{Note.}\\
For how we arrive to this identity, observe that if $q=e^{-\pi\sqrt{r}}$, then $$1+240\sum^{\infty}_{n=1}\frac{n^3q^{2n}}{1-q^{2n}}=u(\beta_r)^4\eqno{(b1)}$$
\[
\]
\textbf{Example.}
$$
u(\beta_{16})=\frac{\sqrt{2} \left(177+124 2^{1/4}+60 \sqrt{2}+68 2^{3/4}\right)^{1/4} \Gamma\left(\frac{5}{4}\right)^2}{\pi^{3/2}}
$$ 
$$\beta_{16}=\frac{16}{761354780+538359129\sqrt{2}+231 \sqrt{3 \left(7242006835334+5120872142664\sqrt{2}\right)}}$$  
\[
\] 
\textbf{Theorem 4.}
\begin{equation}
\beta_r=\frac{1}{2}-\frac{1-20\alpha_{3r}-8\alpha_{3r}^2}{2(1+8 a_{3r})^{3/2}}
\end{equation}
\textbf{Note.}\\
1) Relations (a1),(20),(21),(22),(23) may lead us to the concluding result that exists a number $\xi$ such that $$1-504\sum^{\infty}_{n=1}\frac{n^5q^{2n}}{1-q^{2n}}=(1+\xi \beta_r)u\left(\beta_r\right)^6 .$$  
By setting values one can see that $\xi=-2$.\\
2) Using the triplication formula for $\alpha$ (see [6])
\begin{equation}
\alpha_{9r}=\left(\frac{1-\sqrt[3]{1-\alpha_r}}{1+2\sqrt[3]{1-\alpha_r}}\right)^3
\end{equation}
and in view of Theorem 4 we have if $\alpha'_r=1-\alpha_r$
\begin{equation}
\beta_{3r}=\frac{1}{2}+\frac{27+4\alpha_r(-9+2\alpha_r)}{20(-9+8\alpha_r)(1+2\alpha'^{1/3}_r)^{1/2}\sqrt{1-2\alpha'^{1/3}_r+4\alpha'^{2/3}_r}}
\end{equation}
Also in [6] one can find the 2-degree modular equation for $\alpha_r$ which is $$(\alpha_r\alpha_{4r})^{1/3}+(\alpha'_{r}\alpha'_{4r})^{1/3}=1$$ 

\section{Ramanujan type $1/\pi$ series}

When $q=e^{2\pi i \tau}$, the modular $j$-invariant is defined by
$$j(\tau)=1728\frac{Q^3(q)}{Q^3(q)-R^2(q)}$$
or
$$j(\tau)=\frac{432}{\beta_{r}(1-\beta_{r})}$$
and 
$$t_r=\frac{1}{(1-2\beta_{r/4})u_{r/4}^2}\left(P(q)-\frac{6}{\sqrt{r}\pi}\right)=$$
$$
=\frac{1}{(1-2\beta_{r/4})u_{r/4}^2}\left(3\frac{E(x_{r/4})}{K(x_{r/4})}-2+x_{r/4}-\frac{3\pi}{4\sqrt{r/4}K(x_{r/4})^2}\right)F_{r/4}^2
$$
or
\begin{equation}
t_r=\frac{1+k_{r/4}^2-\frac{6}{\sqrt{r}}a\left(\frac{r}{4}\right)}{\sqrt{1-k_{r/4}^2+k_{r/4}^4}(1-2\beta_{r/4})}
\end{equation}
Also we define 
\begin{equation}
J_r:=4\beta_{r}(1-\beta_{r})
\end{equation} 
\begin{equation}
T_r:=\frac{1+k_{r}^2-\frac{3}{\sqrt{r}}a\left(r\right)}{\sqrt{1-k_{r}^2+k_{r}^4}(1-2\beta_{r})}
\end{equation} 
For the $a(r)$ function see relation (25).
\[
\] 
\textbf{Theorem 5.}(see [7],[8])
\begin{equation}
\frac{3}{\pi\sqrt{r}\sqrt{1-J_{r}}}=\sum^{\infty}_{n=0}\frac{\left(\frac{1}{6}\right)_n\left(\frac{5}{6}\right)_n\left(\frac{1}{2}\right)_n}{(1)^3_n}(J_{r})^n(6n+1-T_r)
\end{equation}
where $(a)_0:=1$ and $(a)_n:=a(a+1)(a+2)\ldots(a+n-1)$, for each positive integer $n$.
\[
\]
\textbf{Examples.}\\
1)For $r=2$
$$
J_2=\frac{27}{125}
$$
$$
T_2=\frac{5}{14}
$$
and
\begin{equation}
\frac{15\sqrt{5}}{14\pi}=\sum^{\infty}_{n=0}\frac{\left(\frac{1}{6}\right)_n\left(\frac{5}{6}\right)_n\left(\frac{1}{2}\right)_n}{(n!)^3}\left(\frac{5}{14}\right)^n\left(6n+\frac{9}{14}\right)
\end{equation}
2) For $r=4$ we have
\begin{equation}
\frac{11 \sqrt{\frac{11}{3}}}{14 \pi }=\sum^{\infty}_{n=0}\frac{\left(\frac{1}{6}\right)_n\left(\frac{5}{6}\right)_n\left(\frac{1}{2}\right)_n}{(n!)^3}\left(\frac{27}{125}\right)^n\left(6n+\frac{10}{21}\right)
\end{equation}
3) For $r=5$ we have
$$T_5=\frac{1}{418} \left(139+45 \sqrt{5}\right)$$
$$J_5=\frac{27 \left(-1975+884 \sqrt{5}\right)}{33275}$$
Hence
$$
\frac{\sqrt{21650+5967 \sqrt{5}}}{\pi }=
$$
\begin{equation}
=\sum^{\infty}_{n=0}\frac{\left(\frac{1}{6}\right)_n\left(\frac{5}{6}\right)_n\left(\frac{1}{2}\right)_n}{(n!)^3}\left(\frac{-53325+23868 \sqrt{5}}{33275}\right)^n\left(836 n+93-15 \sqrt{5}\right) 
\end{equation}
4) For $r=8$ we have 
$$k_8^2=113+80 \sqrt{2}-4 \sqrt{2 \left(799+565 \sqrt{2}\right)}$$
$$a(8)=2 \left(10+7 \sqrt{2}\right) \left(1-\sqrt{-2+2 \sqrt{2}}\right)^2
$$
Then
$$
\frac{15 \sqrt{\frac{5}{2} \left(84125+81432 \sqrt{2}\right)}}{9982 \pi }=
$$
\begin{equation}
=\sum^{\infty}_{n=0}\frac{\left(\frac{1}{6}\right)_n\left(\frac{5}{6}\right)_n\left(\frac{1}{2}\right)_n}{(n!)^3}\left(\frac{5643000-3990168 \sqrt{2}}{1520875}\right)^n\left(\frac{3276-1125 \sqrt{2}+29946 n}{4991}\right)
\end{equation}
5) For $r=18$ we have 
$$k_{18}=(-7+5\sqrt{2})(7-4\sqrt{3})$$
$$a(18)=-3057+2163 \sqrt{2}+1764 \sqrt{3}-1248 \sqrt{6}$$
$$a_{6}=\frac{1}{500}(68-27\sqrt{6})$$
$$\beta_{18}=\frac{1}{2}-\frac{7 \left(49982+4077 \sqrt{6}\right)}{10 \sqrt{5} \left(989+54 \sqrt{6}\right)^{3/2}}$$
$$
J_{18}=\frac{637326171-260186472 \sqrt{6}}{453870144125}
$$
$$T_{18}=\frac{712075+49230 \sqrt{6}}{1074514}$$
6) For $r=27$ $$k_{27}=\frac{1}{2}\sqrt{\frac{1+100\cdot2^{1/3}-80\cdot 2^{2/3}}{2+\sqrt{3-100\cdot2^{1/3}+80\cdot2^{2/3}}}}$$
$$a(27)=3\left[\frac{1}{2}\left(\sqrt{3}+1\right)-2^{1/3}\right]$$
$$
J_{27}=\frac{56143116+157058640\cdot2^{1/3}-160025472\cdot 2^{2/3}}{817400375}
$$
$$T_{27}=\frac{58871825+22512960\cdot2^{1/3}+13208820\cdot 2^{2/3}}{132566687}$$\
\textbf{Some Results.}\\
The duplication formula for $k_r$ is 
$$k_{4r}=\frac{1-k'_r}{1+k'_r}$$
7) For $r=72$ $$k_{72}=\frac{1-\sqrt{-9602+6790 \sqrt{2}+5544 \sqrt{3}-3920 \sqrt{6}}}{1+\sqrt{-9602+6790 \sqrt{2}+5544 \sqrt{3}-3920 \sqrt{6}}}$$
8) For $r=108$ $$k_{108}=\frac{2-\sqrt{2+\sqrt{3-100\cdot 2^{1/3}+80\cdot 2^{2/3}}}}{2+\sqrt{2+\sqrt{3-100\cdot 2^{1/3}+80\cdot 2^{2/3}}}}$$
$$a(108)=\frac{6 \left(1-2\cdot 2^{1/3}-7\cdot \sqrt{3}+4\sqrt{3} \sqrt{2+\sqrt{3-100\cdot 2^{1/3}+80\cdot 2^{2/3}}}\right)}{2+\sqrt{2+\sqrt{3-100\cdot 2^{1/3}+80\cdot 2^{2/3}}}}$$
9) For $r=288$ $$k_{288}=\frac{\left(-1+\left(-9602+6790 \sqrt{2}+5544 \sqrt{3}-3920 \sqrt{6}\right)^{1/4}\right)^2}{\left(1+\left(-9602+6790 \sqrt{2}+5544 \sqrt{3}-3920\sqrt{6}\right)^{1/4}\right)^2}$$
10) For $r=1728$ we set $$p=\sqrt{2+\sqrt{3-100 \cdot2^{1/3}+80\cdot 2^{2/3}}}$$
then
$$k_{108}=\frac{2-p}{2+p}$$ 
$$k_{432}=\left(\frac{\sqrt{p}-\sqrt{2}}{\sqrt{p}+\sqrt{2}}\right)^2\eqno{(k432)}
$$ 
$$
k_{1728}=\left(\frac{2^{3/4} p^{1/4}-\sqrt{2+p}}{2^{3/4} p^{1/4}+\sqrt{2+p}}\right)^2\eqno{(k1728)}
$$
$$
a(108)=\frac{6 \left(4-8\cdot 2^{1/3}+\sqrt{3} p^2\right)}{(2+p)^2}
$$
$$
a(432)=\frac{96-192\cdot 2^{1/3}-48 \sqrt{3}+48 \sqrt{3} p+12 \sqrt{3} p^2}{\left(\sqrt{2}+\sqrt{p}\right)^4}\eqno{(a432)}
$$
$$a(1728)=\frac{384-768\cdot 2^{1/3}-288 \sqrt{3}+192 \sqrt{6} \sqrt{p}-96 \sqrt{3} p+96 \sqrt{6} p^{3/2}+24 \sqrt{3} p^2}{\left(2^{3/4} p^{1/4}+\sqrt{2+p}\right)^4}\eqno{(a1728)}$$
$$\alpha_{144}=\frac{\left(2-\frac{1}{5} 2^{1/6} \left(-3+2 \sqrt{2}+\sqrt{3}+3 \sqrt{6}\right)\right)^3}{8 \left(1+\frac{1}{5} 2^{1/6} \left(-3+2 \sqrt{2}+\sqrt{3}+3 \sqrt{6}\right)\right)^3}\eqno{(\alpha144)}$$
11) From (31) we have
$$  
a_{216}=\frac{\left(46-3\left(53504-36470 \sqrt{2}+29380 \sqrt{3}-20090 \sqrt{6}\right)^{1/3}\right)^3}{8 \left(23+3 \left(53504-36470 \sqrt{2}+29380 \sqrt{3}-20090 \sqrt{6}\right)^{1/3}\right)^3}
$$
Set
$$P=46-3 \left(53504-36470 \sqrt{2}+29380 \sqrt{3}-20090 \sqrt{6}\right)^{1/3}$$
and
$$A=\left(1+\frac{P^3}{8 (-69+P)^3}\right)^{1/3}$$
then
$$
\alpha_{216}=\frac{P^3}{8 (69-P)^3}
$$
and
$$\alpha_{1944}=\left(\frac{1-A}{1+2A}\right)^3$$
Hence
\begin{equation}
J_{5832}=\frac{576 (1-A)^9 A \left(1+A+A^2\right)}{(1+2 A)^3 \left(1+78 A+84 A^2+80 A^3\right)^3}
\end{equation} $$a(81)=\frac{9}{2}\left(1-\sqrt{2}\cdot3^{1/4}(\sqrt{3}+1)(3+y)y^{-1}\right)$$
$$y=\left((2+\sqrt{12})^{1/3}+1\right)^2$$
12)
From the relations ($k432$), ($a432$), ($\alpha144$) we get 
$$J_{432}=\frac{9 (2-w)^9 w (4+2w+w^2)}{64 (1+w)^{3} (1+39 w+21w^2+10 w^3)^3}$$
Where
$$w=\frac{1}{5} 2^{1/6}\left(-3+2\sqrt{2}+\sqrt{3}+3 \sqrt{6}\right)$$
$$p=\sqrt{2+\sqrt{3-100\cdot 2^{1/3}+80\cdot 2^{2/3}}}$$
$$T_{432}=\frac{20-8 \sqrt{3}+16\cdot2^{1/3}\sqrt{3}+12 p+p^2}{\sqrt{16+480 p+536 p^2+120 p^3+p^4} \sqrt{1-J_{432}}}$$
If we set the values $J_{432}$ and $T_{432}$ in Theorem 5 we get a formula that gives 53 digits per term.\\ 
13) It is
$$\alpha_{64}=$$
$$\frac{-6646676+2622114 \sqrt{6}+\sqrt{-76719141489548+31340912395122 \sqrt{6}}}{3375 \left(3375+\sqrt{24683977-5244228 \sqrt{6}-2 \sqrt{-76719141489548+31340912395122 \sqrt{6}}}\right)}$$
Set 
$$
p_1=\frac{1}{2}+\frac{\sqrt{24683977-5244228 \sqrt{6}-2 \sqrt{-76719141489548+31340912395122 \sqrt{6}}}}{6750}
$$
and 
$$p_1=(1-A)^3$$
then from the triplication formula (31) we get  
$$
\alpha_{576}=\left(\frac{A}{3-2A}\right)^3
$$
Hence
$$J_{1728}=\frac{576 (1-A) A^9 \left(3-3 A+A^2\right)}{(2 A-3)^3 \left(-243+486 A-324 A^2+80 A^3\right)^3}
$$
$$T_{1728}=
\frac{(44-16 \sqrt{3}+32 2^{1/3} \sqrt{3}+24 \sqrt{2} \sqrt{p}+36 p+12 \sqrt{2} p^{3/2}-p^2)(1-x)^{-1/2}}{\sqrt{64 p^2+960 \sqrt{2} p^{3/2} (2+p)+1072 p (2+p)^2+120 \sqrt{2} \sqrt{p} (2+p)^3+(2+p)^4}}$$
Setting the values of $T_{1728}$ and $J_{1728}$ in the relation (36) of Theorem 5 we get a formula which gives 110 digits per term.
\newpage
\section{Table of $\alpha_r$}

$$
\alpha_1=\frac{1}{2}
$$
$$
\alpha_2=\frac{1}{4} \left(2-\sqrt{2}\right)
$$
$$
\alpha_3=\frac{1}{4}(3\sqrt{3}-5)
$$
$$
\alpha_4=\frac{1}{18}(9-5\sqrt{3})
$$
$$
\alpha_5=\frac{1}{50}(25-11\sqrt{5})
$$
$$
\alpha_6=\frac{1}{500}(68-27\sqrt{6})
$$
$$
\alpha_7=\frac{-34+13 \sqrt{7}}{6 \left(18+\sqrt{528-78 \sqrt{7}}\right)}
$$
$$
\alpha_8=\frac{1}{265+153 \sqrt{3}+\sqrt{139922+80784 \sqrt{3}}}
$$
$$
\alpha_9=\frac{1}{250} \left(187-171\cdot 2^{1/3}+18 \cdot2^{2/3}\right)
$$
$$
\alpha_{10}=\frac{223-70 \sqrt{10}}{54 \left(54+\sqrt{2470+140 \sqrt{10}}\right)}
$$
$$
\alpha_{11}=\frac{1}{1552+900 \sqrt{3}+10 \sqrt{22 \left(2198+1269 \sqrt{3}\right)}}
$$
$$
\alpha_{12}=\frac{3929-1239 \sqrt{3}-9 \sqrt{-13176+30254 \sqrt{3}}}{5324}
$$
$$
\alpha_{13}=\frac{4}{17743+4921 \sqrt{13}+9 \sqrt{7771398+2155398 \sqrt{13}}}
$$
$$
\alpha_{14}=\frac{1}{7276+1584 \sqrt{21}+10 \sqrt{2 \left(528079+115236 \sqrt{21}\right)}}
$$
$$
\alpha_{15}=\frac{39452-22815 \sqrt{3}+27 \sqrt{5 \left(767653-443204 \sqrt{3}\right)}}{5324}
$$
$$
\alpha_{16}=\frac{2}{37102+15147 \sqrt{6}+45 \sqrt{1359506+555016 \sqrt{6}}}
$$
$$
\alpha_{17}
=(1-115980x+435708x^2-953888x^3+1263024x^4-$$ $$-943296x^5+314432x^6)_1
$$
$$
\alpha_{18}=(1-179016 x+14931960 x^2-21869696 x^3+6249408 x^4-$$
$$
-1116672 x^5+2515456 x^6)_1
$$
$$
\alpha_{19}=\frac{1}{68176+15652 \sqrt{19}+30 \sqrt{6 \left(1722694+395213 \sqrt{19}\right)}}
$$
$$
\alpha_{20}=\frac{2}{205694+118755 \sqrt{3}+\sqrt{84617448935+48853906920 \sqrt{3}}}
$$
$$
\alpha_{21}=\frac{1841863-160683 \sqrt{21}-27 \sqrt{799241022+191433398 \sqrt{21}}}{2456500}
$$
$$
\alpha_{22}=\frac{5636-981 \sqrt{33}-2 \sqrt{15552559-2707353 \sqrt{33}}}{108+6\sqrt{6 -5582+981 \sqrt{33}+2 \sqrt{15552559-2707353 \sqrt{33}}}}
$$
$$
\alpha_{23}=\frac{-9262+5220 \sqrt{3}+9 \sqrt{-342344201+197652864 \sqrt{3}}}{48668+46\sqrt{46}\sqrt{33596-5220 \sqrt{3}-9 \sqrt{-342344201+197652864 \sqrt{3}}}}
$$
$$\alpha_{24}=\frac{-673636+492345 \sqrt{2}-135 \sqrt{6 \left(8352003-5902442 \sqrt{2}\right)}}{48668}$$
$$
\alpha_{25}=\frac{727+436\cdot10^{1/3}-359\cdot 10^{2/3}}{4374+54 \sqrt{3 \left(1460-436\cdot 10^{1/3}+359\cdot 10^{2/3}\right)}}
$$
$$\alpha_{26}=\frac{1}{999652+277200 \sqrt{13}+170 \sqrt{69142562+19176696 \sqrt{13}}}$$
$$
\alpha_{27}=(-1+5686692 x+2583734664 x^2-10597637248 x^3+21994988736 x^4-$$
$$
-22263893760 x^5+8291469824 x^6)_2
$$
$$
\alpha_{28}=-\frac{13954-7631 \sqrt{7}+\sqrt{7 \left(-574473859+219230258 \sqrt{7}\right)}}{39366+162 \sqrt{100911-7631 \sqrt{7}+3 \sqrt{7 \left(-574473859+219230258 \sqrt{7}\right)}}}
$$
$$
\alpha_{29}=(1-11289462 x+11730735 x^2-906935 x^3+514440 x^4-73167 x^5+24389 x^6)_1
$$
$$
\alpha_{30}=(1-15766352 x+12156074560 x^2-76862264960 x^3+187197216640 x^4-$$
$$-234229025792 x^5+166117740544 x^6-66850611200 x^7+12487168000 x^8)_1
$$
$$
\alpha_{31}=\frac{1}{5473576+983164 \sqrt{31}+90 \sqrt{7398145464+1328746146 \sqrt{31}}}
$$
$$
\alpha_{32}=\left(27436813-33403725 \sqrt{3}+171 \sqrt{-46084476791+44878196929 \sqrt{3}}\right)\times $$
$$
\times\left[148035889+12167 \surd \left(93162263+66807450 \sqrt{3}-\right.\right.
$$
$$
\left.\left.-342 \sqrt{-46084476791+44878196929 \sqrt{3}}\right)\right]^{-1}
$$
$$
\alpha_{33}=(1-41580060 x-27923196 x^2-1877533040 x^3-3205978608 x^4-$$
$$
-2128246656 x^5+7668682048 x^6)_1
$$
$$
\alpha_{34}=\frac{1}{14221252+2438800 \sqrt{34}+990 \sqrt{412679650+70773976 \sqrt{34}}}
$$
$$
\alpha_{35}=\left[-765 \sqrt{7 \left(2860957-624312 \sqrt{21}\right)}+2024 \left(-1196+261 \sqrt{21}\right)\right]\times
$$
$$
\times\left[500+10\sqrt{10}\sqrt{2420954-528264 \sqrt{21}+765 \sqrt{7 \left(2860957-624312 \sqrt{21}\right)}}\right]^{-1}
$$
$$
\alpha_{36}=(1-105029898 x+1107006414 x^2-1374218912 x^3+42433836 x^4+$$
$$
+163965000 x^5+166375000 x^6)_1
$$
$$
\alpha_{37}=$$
$$[87483817+8210033\sqrt{37}-9 \sqrt{6 \left(18000682309553+3428999525731 \sqrt{37}\right)}]
$$
$$
[ -118955010-24630099\sqrt{37}+
$$
$$
+27\sqrt{6\left(18000682309553+3428999525731 \sqrt{37}\right)} ]^{-1}
$$
$$
\alpha_{40}=$$
$$\frac{-87049+36036 \sqrt{6}-7 \sqrt{5 \left(61769291-25214724 \sqrt{6}\right)}}{6750+30 \sqrt{1356360-540540 \sqrt{6}+105 \sqrt{5 \left(61769291-25214724 \sqrt{6}\right)}}}
$$
$$
\alpha_{44}=$$
$$\frac{30776-16695\sqrt{3}-117\sqrt{-30041+17490 \sqrt{3}}}{31250+250\sqrt{-15151+16695 \sqrt{3}+117 \sqrt{-30041+17490\sqrt{3}}}}
$$
$$\alpha_{49}=$$
$$\frac{255952117+949411341\cdot 2^{2/3}\cdot7^{1/3}-680830164\cdot 2^{1/3}\cdot 7^{2/3}}{332750 \left(166375+3 \sqrt{3047187612-105490149\cdot 2^{2/3}\cdot 7^{1/3}+75647796\cdot 2^{1/3}\cdot 7^{2/3}}\right)}$$
$$\alpha_{59}=$$
$$\frac{1}{11689799152+6749120700 \sqrt{3}+230 \sqrt{118 \left(43783168572014+25278224160987 \sqrt{3}\right)}}$$
\[
\]
Using the 2,3,5,7 degree modular equations for the cubic base (see [6]) we can evaluate higher order values of $\alpha_r$.

\newpage

\centerline{\bf References}\vskip .2in

[1]: M.Abramowitz and I.A.Stegun. 'Handbook of Mathematical Functions'. Dover Publications

[2]: B.C.Berndt. 'Ramanujan`s Notebooks Part I'. Springer Verlag, New York (1985)

[3]: B.C.Berndt. 'Ramanujan`s Notebooks Part II'. Springer Verlag, New York (1989)

[4]: B.C.Berndt. 'Ramanujan`s Notebooks Part III'. Springer Verlag, New York (1991)

[5]: B.C.Berndt, S.Bhargava and F.G.Garvan. 'Ramanujans theories of elliptic functions to alternative bases'. Trans. Amer. Math. Soc. 347(1995),4163-4244.   

[6]: Bruce C. Berndt and Heng Huat Chan. 'Ramanujan and the Modular j-Invariant'. Canad. Math. Bull. Vol.42(4), 1999. pp.427-440.

[7]: Bruce C. Berndt and Aa Ja Yee. 'Ramanujans Contributions to Eisenstein Series, Especially in his Lost Notebook'. (page stored in the Web). 

[8]: J.M. Borwein and P.B. Borwein. 'Pi and the AGM'. John Wiley and Sons, Inc. New York, Chichester, Brisbane, Toronto, Singapore (1987). 

[9]: I.S. Gradshteyn and I.M. Ryzhik, 'Table of Integrals, Series and Products'. Academic Press (1980).
 
[10]: E.T.Whittaker and G.N.Watson. 'A course on Modern Analysis'. Cambridge U.P. (1927)

[11]: I.J.Zucker. 'The summation of series of hyperbolic functions'. SIAM J. Math. Ana.10.192(1979)

\end{document}